\documentclass{article}

\usepackage{amsmath, amssymb, mathtools}

\title{A geometric decomposition of real diagonalizable matrices with complex eigenvalues.}

\author{Cristobal Arratia}

\date{%
 Nordita\\ 
 KTH Royal Institute of Technology and Stockholm University\\
 Hannes Alfvéns väg 12, SE-106 91 Stockholm, Sweden\\[2ex]
 \today
 }

\begin{document}

\maketitle

\begin{abstract}
For any real diagonalizable matrix with complex eigenvalues we provide a real, coordinate free decomposition with a clear geometric interpretation. 
\end{abstract}

Imaginary numbers have been a source of puzzlement for centuries. 
Square roots of negative numbers had crept into solution formulas for cubic equations in the XVI 
century, and it took some years until Rafael Bombelli managed to show how those formulas lead to real solutions, providing along the way the correct multiplication rules of complex numbers~\cite{CrossleyEmergenceNumber,gonzalez2011journey}. 
This remarkable feat of algebra had to wait for more than two hundred years for its geometric counterpart when, around the turn of the XIX century, multiplication and addition of complex numbers could be visualized in terms of vectors in the complex plane~\cite{gonzalez2011journey,nahin1998imaginary}. 
This includes, of course, the interpretation of the imaginary unit $\sqrt{-1}$ as an operator which rotates the vector of a complex number by $\pi/2$ in the complex plane. 
It would still take until the work of an influential figure like Gauss, who pondered about the use of the term `lateral' instead of `imaginary' or `impossible'~\cite{nahin1998imaginary}, for complex numbers to gain broad acceptance in mainstream Mathematics~\cite{gonzalez2011journey}. 
Despite this acceptance, the question of the side ({\it latus}), if any, to which this `lateral' number relates is often shrugged off. 
More precisely, the visualization of complex numbers in relation to real problems in which they appear is not always clear. 
This is particularly surprising in solutions to problems in which one may expect to be able to think geometrically, like that of the eigenvalues of linear transformations of real vector spaces. 


\subsection*{Matrix Diagonalization}
Consider a matrix $M \in \mathbb{R}^{n\times n}$ with a complete set of eigenvectors $\{b_i\}$ with eigenvalues $\lambda_i$ satisfying
\begin{equation} \label{eq:eigenEq}
  M b_i = \lambda_i b_i, 
\end{equation}
for $i=1,\ldots, n.$ 
The matrix $M$ is diagonalizable~\cite{horn_johnson_2012}, meaning that it can be decomposed as 
\begin{equation} \label{eq:eigenMatrixDiaglze}
  M  = B \Lambda B^{-1}, 
\end{equation}
where $\Lambda$ is a diagonal matrix with the eigenvalues $\lambda_i$ and the columns of $B$ contain the corresponding eigenvectors $b_i.$ 
While the matrix $M$ is real, the eigenvalues and eigenvectors can generally come in complex conjugate pairs. 
Of course, as described in ref. \cite{MurWin31PNAS} in a closely related context nearly a century ago, the decomposition can be written in purely real terms. 
The real canonical form is obtained by replacing in $\Lambda$ the diagonal block of each pair of complex conjugate eigenvalues as
\begin{equation}
 \left(\begin{array}{cc}
        \lambda_i &   0 \\
           0    &   \lambda_i^*
       \end{array} \right) 
       \longrightarrow
      \left(\begin{array}{cc}
         \sigma_i & \omega_i \\
        -\omega_i & \sigma_i
       \end{array} \right),  
\end{equation}
where $\lambda_i=\sigma_i+\sqrt{-1}\,\omega_i,$ and the corresponding $b_i$ and $b_i^*$ in the columns of $B$ by their real and imaginary parts, respectively. 
As it is commonly understood, the imaginary part of the eigenvalue $\omega_i$ has to do with a rotation in the plane determined by the real and imaginary parts of the complex eigenvector $b_i.$ 

A more detailed geometric understanding of the action of the real matrix $M$ on a real vector ${\bf x}$ is attained by considering left eigenvectors $d_i$ satisfying
\begin{equation} \label{eq:leftEigenEq}
 d_i^{\dag} M = \lambda_i d_i^{\dag},
\end{equation}
where $d^\dag=d^{*T}$ is the transpose conjugate of $d.$ 
The left and right eigenvectors satisfy a biorthogonality property whereby $d_i^{\dag}b_j=0$ when $\lambda_i\neq\lambda_j.$ 
The eigenvalue decomposition of the real matrix $M$ can be written as 
\begin{equation}  \label{eq:eigenMatrixAdDecomp}
 M=\sum_{i=1}^n \lambda_i \frac{b_i d_i^{\dag}}{d_i^{\dag} b_i},
\end{equation}
so that its action on a real vector ${\bf x}$ becomes
\begin{equation}
 M{\bf x}= \sum_{i=1}^n \lambda_i b_i \frac{d_i^{\dag}{\bf x}}{d_i^{\dag} b_i}, 
\end{equation}
explicitly showing that the effect of the change of base matrix $B^{-1}$ in \eqref{eq:eigenMatrixDiaglze} is the scalar projection of $\bf x$ along the left eigenvectors normalized by the scalar product between the right and left eigenvectors $d_i^{\dag} b_i.$ 
For the case in which the eigenvalues are complex the geometric interpretation of these straightforward operations is not completely transparent. 

\section*{Geometric $\sqrt{-1}.$}

A geometric understanding of possible origins of a $\sqrt{-1}$ in an $n-$dimensional real vector space is provided by Geometric Algebra. 
Extensively developed and forcefully advocated by Hestenes (see for example \cite{HestenesReform} and references therein), Geometric Algebra has attracted efforts from a variety of scholars who have embarked in its development and promotion in different contexts such as Physics~\cite{gull1993imaginary,DoranLasenby2003}, Computer Science~\cite{Dorst2009geometric} and undergraduate education~\cite{macdonaldTextbook}. 
While this is no place for an exposition of Geometric Algebra (see \cite{macdonald2017survey} for a brief survey), we give a brief description of some of its aspects which are relevant here. 
Geometric Algebra (GA) is an associative algebra (a Clifford algebra) whose product is the geometric product, which in GA is implied by yuxtaposition but we will denote it here by $\odot.$ 
Its result between two vectors yields the (direct) sum of a scalar and a {\it 2-blade} or {\it bivector} as  
\begin{equation} \label{eq:geomProduct}
 {\bf u}\odot{\bf v} = {\bf u}\cdot{\bf v} + {\bf u} \wedge {\bf v},
\end{equation}
where $\cdot$ is the usual dot product and $\wedge$ represents the exterior or wedge product, commonly referred to in GA as the outer product (not to be confused with the tensor product). 
For our purposes, it will suffice to define it in terms of standard matrix multiplication as
\begin{equation} \label{eq:wedgeDef}
 {\bf u} \wedge {\bf v} \equiv {\bf uv}^T-{\bf vu}^T.
\end{equation}
A 2-blade or bivector can be geometrically interpreted as an oriented area, that is, a piece of the plane spanned by ${\bf u}$ and ${\bf v}$ with an associated number (with magnitude equal to the area of the parallelogram defined by ${\bf u}$ and ${\bf v}),$ but with no particular shape. 
In addition, a 2-blade can also represent the plane $\operatorname{span}({\bf u},{\bf v});$ it is a (real) scalar factor of the outer product of any pair of vectors forming a base of that plane. 

As exemplified after eq. \ref{eq:wedgeDef}, the exterior product $\wedge$ generates $k-$blades, geometrical objects representing $k-$dimensional subspaces $(k\leq n)$ spanned by its linearly independent factors. 
Here we will not need to go beyond $k=2.$  
As exemplified by eq. \ref{eq:geomProduct}, the geometric product $\odot$ generates different types of geometrical objects from its factors and adds them together; the example in eq. \ref{eq:geomProduct} includes the addition of a scalar projection and a $2-$blade. 
While the addition of two quantitites of different nature may appear unusual, it is similar to the addition of real and imaginary parts of complex numbers. 
In fact, the similiarities between GA and complex numbers go further. 
If ${\bf e}_1$ and ${\bf e}_2$ are any pair of orthonormal vectors in $\mathbb{R}^n,$ then
\begin{align}
  ({\bf e}_1\odot{\bf e}_2)\odot({\bf e}_1\odot{\bf e}_2)&=({\bf e}_1\odot{\bf e}_2)\odot(-{\bf e}_2\odot{\bf e}_1) \nonumber\\
 &=-{\bf e}_1\odot({\bf e}_2\odot{\bf e}_2)\odot{\bf e}_1\\
 &=-({\bf e}_2\cdot{\bf e}_2){\bf e}_1\odot{\bf e}_1=-{\bf e}_1\cdot{\bf e}_1=-1, \nonumber
\end{align}
where we have used eqs. \ref{eq:geomProduct} and \ref{eq:wedgeDef}, and the associativity of the geometric product. 
Therefore, under the geometric product, any 2-blade in GA is proportional to a $\sqrt{-1}.$ 
The fact that GA gives a geometric origin and interpretation to the imaginary unit is often highlighted by advocates of GA; the geometric products \ref{eq:geomProduct} of any two vectors in a given plane have been called `GA complex numbers' since they are isomorphic to the usual complex numbers~\cite{macdonald2017survey}. 
This, among other possible types of $\sqrt{-1}$ in GA, is the kind of imaginary unit that we will consider. 
These 2-blades, together with the usual scalars and vectors, are the only objects from GA that we will need; 
there will be no further explicit use of the geometric product $\odot.$ 

For any vector ${\bf x},$ and any pair of 
vectors ${\bf u}$ and ${\bf v},$ 
we have 
\begin{equation} \label{eq:bladeAction}
 {\bf x}\cdot[({\bf u}\wedge{\bf v}) {\bf x}]=0. 
\end{equation}
We mention two aspects about the action of ${\bf u}\wedge{\bf v}$ on ${\bf x}$ that appears in the square bracket. 
First, as seen from definition \ref{eq:wedgeDef}, it is a vector contained in the plane $\operatorname{span}({\bf u},{\bf v}),$ with no contribution from ${\bf x}_{\perp},$ the part of ${\bf x}$ which is orthogonal to that plane. 
Second, ${\bf x}_{\parallel},$ the part of ${\bf x}$ in the plane $\operatorname{span}({\bf u},{\bf v}),$ is made orthogonal to itself by the multiplication by ${\bf u}\wedge{\bf v};$ therefore the blade rotates ${\bf x}_{\parallel}$ by $\pi/2.$ 
This is the role of the imaginary unit as an operator on complex vectors in the complex plane, but it happens here with real vectors on a real geometric plane. 
Inspection shows that this rotation is in the direction that brings ${\bf v}$ towards ${\bf u}.$

Given that much emphasis to geometric interpretation of imaginary units is given in GA, one expects the subject of complex eigenvalues to have been addressed. 
This has been done by extending the notion of eigenvector to that of {\it eigenblade}, which represent invariant subspaces associated to linear operators (whose domain has been extended to all elements of GA, \cite{hestenes1984clifford,DoranLasenby2003}). 
The eigenvalues associated to such eigenblades are real, and the precise geometric connection to the usual complex eigenvalues and eigenvectors has, to the best of our knowledge, not been explicitly established. 
This is possibly the case because GA does not directly build simple linear operators like ${\bf u}{\bf v}^T,$ which involve the symmetric part of the tensor product. 
While we borrow from GA the exterior (or outer) product as defined in eq. \ref{eq:wedgeDef} and its interpretation, we proceed with standard matrix notation to make the result more broadly accessible. 

%

\section*{Real Decomposition with Complex Eigenvalues} 

Let us consider a pair of complex conjugate eigenvalues $(\lambda,\lambda^{*})$ with their corresponding right and left eigenvectors $(b,b^{*})$ and $(d,d^{*})$ satisfying eqs. \ref{eq:eigenEq} and \ref{eq:leftEigenEq}, respectively. 
We write the right and left eigenvectors in terms of their real and imaginary parts as
\begin{subequations}
\begin{align}
 b &= {\bf b} + \sqrt{-1}\, {\bf p},\\
 d &= {\bf d} + \sqrt{-1}\, {\bf q}.
\end{align}
\end{subequations}
The orthogonality satisfied between $d^*$ and $b$, i.e., $d^{*\dag}b=0,$ translates into the conditions 
\begin{subequations}\label{eq:biorthogonalities}
\begin{align}
 {\bf d} \cdot {\bf b} &= {\bf q} \cdot {\bf p}, \\
 {\bf q} \cdot {\bf b} &= -{\bf d}\cdot {\bf p}.
\end{align}
\end{subequations}
If $b$ 
and $d$ 
were associated to different pairs of complex conjugate eigenvalues, the biorthogonality property 
implies that each side of these equations would be zero. 

The four conditions that define the right (eq. \ref{eq:eigenEq}) and left (eq. \ref{eq:leftEigenEq}) eigenvectors associated to the same complex eigenvalue $\lambda=\sigma + \sqrt{-1}\,\omega$ are simultaneously satisfied by
\begin{subequations}\label{eq:realDecomposition}
 \begin{align}
 M_{\lambda}&\equiv \frac{1}{\nu}({\bf b} \wedge {\bf p})\left[\omega \{{\bf dd}^T+{\bf qq}^T\} - \sigma ({\bf d} \wedge {\bf q})\right],\label{eq:realDecLeftBlade}\\
  &= \frac{1}{\nu}\left[\omega \{{\bf bb}^T+{\bf pp}^T\} - \sigma ({\bf b} \wedge {\bf p})\right]({\bf d} \wedge {\bf q}),\label{eq:realDecRightBlade}
\end{align}
\end{subequations}
where the normalization is given by 
\begin{equation}
 \nu\equiv({\bf d\cdot b})^2+ ({\bf q\cdot b})^2.
\end{equation}
The equality between the right hand side expressions in eqs. \ref{eq:realDecLeftBlade} and \ref{eq:realDecRightBlade} follows from conditions \ref{eq:biorthogonalities}. 
We note that every bracket in \ref{eq:realDecomposition} is invariant to the complex phase of the eigenvectors, and that the same expression is obtained if written in terms of the complex conjugate eigenvalue $\lambda^*.$ 
Interchanging the positions of the vectors $({\bf b},{\bf p})$ of the right eigenspace with those of the vectors $({\bf d},{\bf q})$ of the left corresponds to a transposition $(\,\cdot\,)^T$ and the substitution $\omega\rightarrow-\omega.$

There is a clear geometric interpretation of eqs. \ref{eq:realDecomposition} for the action of $M_{\lambda}$ on real vectors. 
To aid the interpretation, it is useful to keep in the back of one's mind the case of a normal matrix (i.e. when $MM^{T}=M^{T}M),$ for which one can identify the right and left eigenvectors. 
As seen from eq. \ref{eq:realDecRightBlade}, a column vector multiplied on the right is first projected on the plane of the left eigenvector represented by $({\bf d} \wedge {\bf q})$ and rotated by $\pi/2$ in the direction going from ${\bf q}$ to ${\bf d}.$ 
The resulting vector is then brought to the plane of the right eigenvector in two different ways, as given by the two terms in the square bracket of eq. \ref{eq:realDecRightBlade}. 
The real part $\sigma$ of the eigenvalue multiplies the vector as is rotated back towards its original orientation by $-({\bf b} \wedge {\bf p}).$  
The imaginary part $\omega$ of the eigenvalue multiplies the term in the curly brackets which is similar to a projection on the plane $\operatorname{span}({\bf b},{\bf p}),$ thus keeping the vector close to its rotated orientation. 
In the case of a normal matrix, ${\bf b}$ and ${\bf p}$ are orthogonal to each other and of the same norm (as seen from conditions \ref{eq:biorthogonalities} with $d=b),$ so the curly brackets in eq. \ref{eq:realDecomposition} are proportional to the identity of the plane; the orientation is then preserved exactly after the first rotation by $\pi/2.$ 
As seen in eq. \ref{eq:realDecLeftBlade}, essentially the same description applies to a multiplication by a row vector on the left. 

As mentioned above, $M_\lambda$ given in eq. \ref{eq:realDecomposition} reproduces the action of the full matrix $M$ on (left and right) eigenvectors associated to a complex eigenvalue $\lambda$ (eqs. \ref{eq:eigenEq} and \ref{eq:leftEigenEq}). 
The biorthogonality property guarantees that right (left) multiplication of $M_{\lambda}$ by other right (resp. left) eigenvectors of $M$ is zero (which may require orthogonalization when there are repeated eigenvalues). 
If the basis of eigenvectors is complete, then eq. \ref{eq:realDecomposition} can be used to provide a real geometric decomposition of the full matrix $M$ by replacing the complex conjugate pairs in eq. \ref{eq:eigenMatrixAdDecomp}. 
To be explicit, and using ${\bf a}_i$ and ${\bf c}_i$ to denote right and left eigenvectors associated to, say, $j$ real eigenvalues $\alpha_i,$ a diagonalizable matrix can be decomposed as
\begin{equation}
 M=\sum_{i=1}^j \alpha_i \frac{{\bf a}_i {\bf c}_i^T}{{\bf c}_i\cdot {\bf a}_i}+ \sum_{i=1}^k M_{\lambda_i},
\end{equation}
where $k$ is the number of complex conjugate pairs so that $j+2k=n,$ and the $M_{\lambda_i}$ are given by eq. \ref{eq:realDecomposition}. 
This provides a real, coordinate free representation of real diagonalizable linear operators. 

\section*{Acknowledgments}
The author acknowledges support of Nordita and the Swedish Research Council Grant No. 2018-04290; Nordita is partially supported by Nordforsk.


\end{document}